\documentclass[12pt]{article}
\usepackage{mathrsfs}
\usepackage{amssymb,a4}
\usepackage{amsmath,amsfonts,amssymb,amsthm}
\newcommand{\h}{\frak{h}}
\newcommand{\al}{\alpha}

\newcommand{\1}{{{\bf 1}}}

\newcommand{\End}{{\rm End}\,}

\def\mod{{\rm mod\,}}

\def\wt{{\rm wt}}

\def\be{\beta}

\def\C{{\mathbb C}}
\def\Q{{\mathbb Q}}

\def\Z{{\mathbb Z}}

\def\1{{\bf 1}}

\def \wt{{\rm wt}}

\def \End{{\rm End}}

\def \mod{{\rm mod}}

\def \<{\langle}
\def \>{\rangle}
\def \w{\omega}
\def \pf{\noindent {\bf Proof: \,}}
\def\theequation{5.\arabic{equation}}

\renewcommand{\theequation}{\thesection.\arabic{equation}}

\newtheorem{theorem}{Theorem}[section]
\newtheorem{prop}[theorem]{Proposition}
\newtheorem{lem}[theorem]{Lemma}

\theoremstyle{definition}
\newtheorem{definition}[theorem]{Definition}

\begin{document}
\begin{center}
{\Large {\bf Representations of the vertex operator
algebra $V_{L_{2}}^{A_{4}}$}} \\

\vspace{0.5cm} Chongying Dong\footnote{Supported by NSF grants and a
Faculty research grant from  the University of California at Santa
Cruz.}
\\
Department of Mathematics,  University of California, Santa Cruz, CA 95064 \\
\vspace{.1 cm} Cuipo Jiang\footnote{Supported  by China NSF grants
10931006, the RFDP of China(20100073110052), and the Innovation
Program of Shanghai Municipal
Education Commission (11ZZ18).}\\
 Department of Mathematics, Shanghai Jiaotong University, Shanghai 200240 China
\end{center}
\hspace{1cm}
\begin{abstract} The rationality and $C_2$-cofiniteness of the orbifold vertex operator algebra $V_{L_{2}}^{A_{4}}$ are
 established and all the irreducible modules are constructed and classified. This is part of classification of rational vertex operator algebras with $c=1.$

2000MSC:17B69
\end{abstract}

\section{Introduction}

Motivated by the classification of rational vertex operator algebras
with $c=1,$ we investigate the vertex operator algebra
$V_{L_{2}}^{A_{4}}$ where $L_2$ is the root lattice of type $A_1$
and $A_4$ is the alternating group  which is a subgroup of the
automorphism group of lattice vertex operator algebra $V_{L_2}.$ The
$C_2$-cofiniteness and rationality of $V_{L_{2}}^{A_{4}}$ are
obtained, and the irreducible modules are constructed and
classified.

Classification of rational vertex operator algebras with $c=1$ goes
back to \cite{G} and \cite{K} in the literature of physics at
character level under the assumption that each irreducible character
is a modular function over a congruence subgroup and the sum of the
square norm of irreducible characters is invariant under the modular
group. According to \cite{K}, the character of a rational vertex
operator algebra with $c=1$ must be the character of one of the
following vertex operator algebras: (a) lattice vertex operator
algebras $V_L$ associated
 with positive definite even lattices $L$ of rank one, (b) orbifold vertex operator algebras $V_L^+$ under the automorphism of $V_L$ induced
 from the $-1$ isometry of $L$, (c) $V_{\Z\alpha}^G$
where $(\alpha,\alpha)=2$ and $G$ is a finite subgroup of $SO(3)$
isomorphic to one of $\{A_4, S_4, A_5\}.$  As it is pointed out in
\cite{DJ1} that this list is not correct if the effective central
charge $\tilde{c}$ \cite{DM2} is not equal to $c.$ The vertex
operator algebra $V_L$ for any positive definite even lattice $L$
has been characterized by using $c,$ the effective central charge
$\tilde{c}$ and the rank of the weigh one subspace as a Lie algebra
\cite{DM2}.
 The orbifold vertex operator algebras $V_L^+$ for rank one lattices $L$ have also been characterized in \cite{DJ1}-\cite{DJ3} and \cite{ZD}.
  But the vertex operator algebra
  $V_{\Z\alpha}^G$ has not been understood well as $G$ is not a cyclic group. Although $V_{\Z\alpha}^G$ is in the above list of rational vertex
operator algebras, the rationality of $V_{\Z\alpha}^G$ was unknown.
The present paper deals with the case $G=A_4.$

The main idea is to realize $V_{\Z\alpha}^G$ as
$(V_{\Z\beta}^+)^{\<\sigma\>}$ where $(\beta,\beta)=8$ and $\sigma$
is an automorphism of
 $sl(2,\C)$ of order $3.$ The vertex operator algebra $V_{\Z\beta}^+$ is well understood (see \cite{DN1}-\cite{DN3}, \cite{A1}-\cite{A2}). Also it is
 easier to deal with the cyclic group $\langle \sigma\rangle$ than nonabelian group $A_4.$ One key step is to give an explicit expression of
  the generator $u^{(9)}$ of weight $9.$  Another key step is to
  prove the
   $C_2$-cofiniteness of $V_{\Z\al}^{A_{4}}$. We achieve this  by using the fusion rules of the Virasoro vertex operator algebra
 $L(1,0)$ and technical calculations. The rationality follows
   from the $C_2$-cofiniteness \cite{M2}. For the classification  of irreducible modules, we follow the standard procedure. We first construct the
   irreducible $\sigma^i$-twisted $V_{\Z\beta}^+$-modules, and
then give the irreducible $(V_{\Z\beta}^+)^{\<\sigma\>}$-submodules.
According to \cite{M1}, these irreducible modules should give a
complete list of irreducible $(V_{\Z\beta}^+)^{\<\sigma\>}$-modules.

It is expected that  the ideas and techniques developed in this
paper will work for  $V_{L_2}^{S_{4}}$  as well. The case $G=A_{5}$
might be more complicated. Once the rationality of $V_{L_2}^{G}$ is
established for all $G,$ the classification of rational vertex
operator algebras with $c=1$ is equivalent to  the following
conjecture: If $V$ is a simple, rational vertex operator algebra of
CFT type such that $\dim V_4<3$ then $V$ is isomorphic to
$V_{L_2}^{G}$ for $G=A_4, S_4, A_5.$

The paper is organized as follows. We recall various notions of
twisted modules from \cite{DLM1} in Section 2. We also briefly
discuss lattice vertex operator algebras $V_L$ \cite{FLM} and
$V_L^+$ including the classification of irreducible modules and
rationality \cite{DN1}-\cite{DN3}, \cite{A2}, \cite{AD}, \cite{DJL}.
In Section 3, we identify the vertex operator algebra
$V_{L_2}^{A_4}$ with $(V_{\Z\beta}^+)^{\<\sigma\>}$ and discuss
several special vectors (which play important roles in later
sections) in both $V_{\Z\beta}^+$ and
$(V_{\Z\beta}^+)^{\<\sigma\>}.$ The rationality and
$C_2$-cofiniteness of $(V_{\Z\beta}^+)^{\<\sigma\>}$ are established
in Section 4. The classification of the irreducible
$(V_{\Z\beta}^+)^{\<\sigma\>}$-modules is achieved in Section 5.

\section{Preliminaries}
\def\theequation{2.\arabic{equation}}
\setcounter{equation}{0}

We first recall weak twisted-modules and twisted-modules for vertex
operator algebras from \cite{DLM2}. Let $(V,Y,{\bold 1},\omega)$ be
a vertex operator algebra \cite{B}, \cite{FLM} and $g$  an
automorphism of $V$ of finite order $T.$ Denote the decomposition of
$V$ into eigenspaces with respect to the action of $g$ as
\begin{equation}\label{g2.1}
V=\oplus_{r\in \Z/T\Z}V^r
\end{equation}
where $V^r=\{v\in V|gv=e^{-2\pi ir/T}v\}$.

\begin{definition} A {\em weak $g$-twisted $V$-module} $M$ is a vector space equipped
with a linear map
$$\begin{array}{l}
V\to (\End\,M)\{z\}\\
v\mapsto\displaystyle{ Y_M(v,z)=\sum_{n\in\Q}v_nz^{-n-1}\ \ \
(v_n\in \End\,M)}
\end{array}$$
which satisfies the following for all $0\leq r\leq T-1,$ $u\in V^r$,
$v\in V,$ $w\in M$,
\begin{eqnarray}
& &Y_M(u,z)=\sum_{n\in \frac{r}{T}+\Z}u_nz^{-n-1} \label{1/2}\\
& &u_lw=0\ \ \
\mbox{for}\ \ \ l>>0\label{vlw0}\\
& &Y_M({\bold 1},z)=1;\label{vacuum}
\end{eqnarray}
 \begin{equation}\label{jacobi}
\begin{array}{c}
\displaystyle{z^{-1}_0\delta\left(\frac{z_1-z_2}{z_0}\right)
Y_M(u,z_1)Y_M(v,z_2)-z^{-1}_0\delta\left(\frac{z_2-z_1}{-z_0}\right)
Y_M(v,z_2)Y_M(u,z_1)}\\
\displaystyle{=z_1^{-1}\left(\frac{z_2+z_0}{z_1}\right)^{r/T}
\delta\left(\frac{z_2+z_0}{z_1}\right) Y_M(Y(u,z_0)v,z_2)}.
\end{array}
\end{equation}
\end{definition}
It is known that (see \cite{DLM2}, etc) the twisted-Jacobi identity
is equivalent to the following two identities.
$$
[u_{m+{r\over T}}, v_{n+{s\over T}}]= \sum_{i=0}^{\infty}{m+{r\over
T}\choose i}(u_iv)_{m+n+{r+s\over T}-i},
$$
$$
\sum_{i\geq 0}{s\choose i}(u_{m+i}v)_{n+\frac{s+t}{T}-i}=\sum_{i\geq
0}(-1)^i {m\choose
i}(u_{m+\frac{s}{T}-i}v_{n+\frac{t}{T}+i}-(-1)^{m}v_{m+n+\frac{t}{T}-i}u_{\frac{s}{T}+i}),$$
where $p,n\in\Z$, $u\in V^s, v\in V^t$.
\begin{definition}
An \textit{admissible $g$-twisted $V$-module} $M=\bigoplus_{n\in
\frac{1}{T}\Z_+} M(n)$ is a $\frac{1}{T}\Z_+$-graded weak
$g$-twisted module such that $u_mM(n)\subset M(\wt u-m-1+n)$ for
$u\in V$ and $m,n\in\frac{1}{T}\Z.$
\end{definition}

\begin{definition}
A \textit{$g$-twisted $V$-module} $M=\bigoplus_{\lambda\in\C}
M_{\lambda}$ is a $\C$-graded weak $g$-twisted $V$-module with
$M_\lambda=\{u\in M|L(0)u=\lambda u\}$ such that $M_\lambda$ is
finite dimensional and for fixed $\lambda\in\C$, $M_{\lambda+n/T}=0$
for sufficiently small integer $n$.
\end{definition}

We now review the vertex operator algebras $M(1)^+,$ $V_L^+$ and
related results from \cite{A1}, \cite{A2}, \cite{AD}, \cite{ADL},
\cite{DN1}, \cite{DN2}, \cite{DN3}, \cite{DJL},
 \cite{FLM}.

Let $L=\Z \alpha$ be a positive definite even lattice of rank one.
That is, $(\alpha,\alpha)=2k$ for some positive integer $k.$  Set
$\h=\C\otimes_{\Z} L$ and extend $(\cdot\,,\cdot)$ to a
$\C$-bilinear form on $\h$. Let
$\hat{\h}=\C[t,t^{-1}]\otimes\h\oplus\C K$ be the affine Lie algebra
associated to the abelian Lie algebra $\h$ so that
\begin{align*}
[\alpha(m),\,\alpha(n)]=2km\delta_{m+n,0}K\hbox{ and }[K,\hat{\h}]=0
\end{align*}
for any $m,\,n\in\Z$, where $\alpha(m)=\alpha\otimes t^m.$  Then
$\hat{\h}^{\geq 0}=\C[t]\otimes\h\oplus\C K$ is a commutative
subalgebra. For any $\lambda\in\h$, we  define a one-dimensional
$\hat{\h}^{\geq 0}$-module $\C e^\lambda$ such that $\alpha(m)\cdot
e^\lambda=(\lambda,\alpha)\delta_{m,0}e^\lambda$ and $K\cdot
e^\lambda=e^\lambda$ for $m\geq0$. We denote by
\begin{align*}
M(1,{\lambda})=U(\hat{\h})\otimes_{U(\hat{\h}^{\geq 0})}\C
e^\lambda\cong S(t^{-1}\C[t^{-1}])\ ({\rm linearly})
\end{align*}
the $\hat{\h}$-module induced from $\hat{\h}^{\geq 0}$-module $\C
e^\lambda$. Set
$$M(1)=M(1,0).$$ Then there exists a linear map
$Y:M(1)\to\End M(1)[[z,z^{-1}]]$ such that $(M(1),\,Y,\,\1,\,\w)$
carries a simple vertex operator algebra structure and
$M(1,\lambda)$ becomes an irreducible $M(1)$-module for any
$\lambda\in\h$ (see \cite{FLM}). The vacuum vector and the Virasoro
element are given by $\1=e^0$ and $\w=\frac{1}{4k}\alpha(-1)^2\1,$
respectively.

Let $\C[L]$ be the group algebra of $L$ with a basis $e^{\beta}$ for
$\beta\in L.$ The lattice vertex operator algebra associated to $L$
is given by
$$V_L=M(1)\otimes \C[L].$$
The dual lattice $L^{\circ}$ of $L$ is
$$L^{\circ}=\{\,\lambda\in\h\,|\,(\alpha,\lambda)\in\Z\,\}=\frac{1}{2k}L.$$
Then $L^{\circ}=\cup_{i=-k+1}^k(L+\lambda_i)$ is the coset
decomposition with $\lambda_i=\frac{i}{2k}\alpha.$ In particular,
$\lambda_0=0.$ Set $\C[L+\lambda_i]=\bigoplus_{\beta\in L}\C
e^{\beta+\lambda_i}.$ Then each $\C[L+\lambda_i]$ is an
$L$-submodule in an obvious way. Set
$V_{L+\lambda_i}=M(1)\otimes\C[L+\lambda_i]$. Then $V_L$ is a
rational vertex operator algebra and $V_{L+\lambda_i}$ for
$i=-k+1,\cdots,k$ are the irreducible modules for $V_L$ (see
\cite{B}, \cite{FLM}, \cite{D1}).

Define a linear isomorphism $\theta:V_{L+\lambda_i}\to
V_{L-\lambda_i}$ for $i\in\{-k+1,\cdots,k\}$ by
\begin{align*}
\theta(\alpha(-n_{1})\alpha(-n_{2})\cdots \alpha(-n_{s})\otimes
e^{\beta+\lambda_i})=(-1)^{k}\alpha(-n_{1})\alpha(-n_{2})\cdots
\alpha(-n_{s})\otimes e^{-\beta-\lambda_i}
\end{align*}
where $n_j>0$ and $\beta\in L.$  Then $\theta$ defines a linear
isomorphism from $V_{L^{\circ}}=M(1)\otimes \C[L^{\circ}]$ to itself
such that
 $$\theta(Y(u,z)v)=Y(\theta u,z)\theta v$$
for $u\in V_{L}$ and $v\in V_{L^{\circ}}.$ In particular, $\theta$
is an automorphism of $V_{L}$ which induces an automorphism of
$M(1).$

For any $\theta$-stable subspace $U$ of $V_{L^{\circ}}$, let $U^\pm$
be the $\pm1$-eigenspace of $U$ for $\theta$. Then $V_L^+$ is a
simple vertex operator algebra.

Also recall the $\theta$-twisted Heisenberg algebra $\h[-1]$ and its
irreducible module $M(1)(\theta)$ from \cite{FLM}. Let $\chi_s$ be a
character of $L/2L$ such that $\chi_s(\alpha)=(-1)^s$ for $s=0,1$
and $T_{\chi_s}=\C$ the irreducible $L/2L$-module with character
$\chi_s$. It is well known that
$V_L^{T_{\chi_s}}=M(1)(\theta)\otimes T_{\chi_s}$ is an irreducible
$\theta$-twisted $V_L$-module (see \cite{FLM}, \cite{D2}). We define
actions of $\theta$ on  $M(1)(\theta)$ and $V_L^{T_{\chi_s}}$ by
\begin{align*}
\theta(\alpha(-n_{1})\alpha(-n_{2})\cdots
\alpha(-n_{p}))=(-1)^{p}\alpha(-n_{1})\alpha(-n_{2})\cdots
\alpha(-n_{p})
\end{align*}
\begin{align*}
\theta(\alpha(-n_{1})\alpha(-n_{2})\cdots \alpha(-n_{p})\otimes
t)=(-1)^{p}\alpha(-n_{1})\alpha(-n_{2})\cdots \alpha(-n_{p})\otimes
t
\end{align*}
for $n_j\in \frac{1}{2}+\Z_{+}$ and $t\in T_{\chi_s}$. We denote the
$\pm 1$-eigenspaces of $M(1)(\theta)$ and $V_L^{T_{\chi_s}}$ under
$\theta$ by $M(1)(\theta)^{\pm}$ and $(V_L^{T_{\chi_s}})^{\pm}$
respectively. We have the following results:
\begin{theorem}\label{t32}
Any irreducible module for the vertex operator algebra $M(1)^+$ is
isomorphic to one of the following modules:$$ M(1)^+, M(1)^-, M(1,
\lambda) \cong M(1, -\lambda)\ (0\neq \lambda \in \h),
M(1)(\theta)^+, M(1)(\theta)^- .$$
\end{theorem}
\begin{theorem}\label{t33}
Any irreducible $V_L^+$-module is isomorphic to one of the following
modules:
$$V_L^{\pm}, V_{\lambda_i+L}( i \not= k),
V_{\lambda_k+L}^{\pm}, (V_L^{T_{\chi_s}})^{\pm}.$$
\end{theorem}

\begin{theorem}\label{t34}  $V_{L}^{+}$ is rational.
\end{theorem}

We remark that the classification of irreducible modules for
arbitrary $M(1)^+$ and $V_L^+$ are obtained in \cite{DN1}-\cite{DN3}
and \cite{AD}. The rationality of $V_L^+$ is established in
\cite{A2} for rank one lattice  and \cite{DJL} in general.

\vskip 0.3cm We next turn our attention to the fusion rules of
vertex operator algebras. Let $V$ be a vertex operator algebra, and
$ W^i$ $ (i=1,2,3$)
 be  ordinary $V$-modules. We denote by $I_{V} \left(\hspace{-3 pt}\begin{array}{c} W^3\\
W^1\,W^2\end{array}\hspace{-3 pt}\right)$  the vector space of all
intertwining operators of type $\left(\hspace{-3 pt}\begin{array}{c}
W^3\\ W^1\,W^2\end{array}\hspace{-3 pt}\right)$.
 For a $V$-module $W$, let
$W^{\prime}$ denote the graded dual of $W$. Then $W'$ is also a
$V$-module \cite{FHL}. It is well known that fusion rules have the
following symmetry (see \cite{FHL}).

\begin{prop}\label{p4.2}
Let $W^{i}$ $(i=1,2,3)$ be $V$-modules. Then
$$\dim I_{{V}} \left(\hspace{-3 pt}\begin{array}{c} W^3\\
W^1\,W^2\end{array}\hspace{-3 pt}\right)=\dim I_{{V}} \left(\hspace{-3 pt}\begin{array}{c} W^3\\
W^2\,W^1\end{array}\hspace{-3 pt}\right), \ \ \ \dim I_{{V}} \left(\hspace{-3 pt}\begin{array}{c} W^3\\
W^1\,W^2\end{array}\hspace{-3 pt}\right)=\dim I_{{V}} \left(\hspace{-3 pt}\begin{array}{c} (W^2)^{\prime}\\
W^1\,(W^3)^{\prime}\end{array}\hspace{-3 pt}\right).$$
\end{prop}

Recall that $L(c,h)$ is the irreducible highest weight module for
the Virasoro algebra with central charge $c$ and highest weight $h$
for $c,h\in \C.$ It is well known that $L(c,0)$ is a vertex operator
algebra. The following two results were obtained in \cite{M} and
\cite{DJ1}.
\begin{theorem}\label{2t1} (1) We have
$$
\dim I_{L(1,0)} \left(\hspace{-3 pt}\begin{array}{c} L(1,k^{2})\\
L(1, m^{2})\,L(1, n^{2})\end{array}\hspace{-3 pt}\right)=1,\ \
k\in{\mathbb Z}_{+},  \ |n-m|\leq k\leq n+m,$$
$$\dim I_{L(1,0)} \left(\hspace{-3 pt}\begin{array}{c} L(1,k^{2})\\
L(1, m^{2})\,L(1, n^{2})\end{array}\hspace{-3 pt}\right)=0,\ \
k\in{\mathbb Z}_{+},  \ k<|n-m| \ {\rm or} \  k>n+m, $$ where
$n,m\in{\mathbb Z}_{+}$.

(2) For $n\in{\mathbb Z}_{+}$ such that $n\neq p^{2}$, for all
$p\in{\mathbb Z}_{+}$, we have
$$\dim I_{L(1,0)} \left(\hspace{-3 pt}\begin{array}{c} L(1,n)\\
L(1, m^{2})\,L(1, n)\end{array}\hspace{-3 pt}\right)=1,$$
$$\dim I_{L(1,0)} \left(\hspace{-3 pt}\begin{array}{c} L(1,k)\\
L(1, m^{2})\,L(1, n)\end{array}\hspace{-3 pt}\right)=0,$$
 for $k\in{\mathbb Z}_{+}$ such that  $k\neq n$.

\end{theorem}

\section{The vertex operator subalgebra $V_{L_{2}}^{A_{4}}$}
\def\theequation{3.\arabic{equation}}
\setcounter{equation}{0} Let $L_{2}=\Z\al$ be the rank one
positive-definite even lattice such that $(\al,\al)=2$ and
$V_{L_{2}}$ the associated simple rational vertex operator algebra.
Then $(V_{L_{2}})_{1}\cong sl_{2}(\C)$ and $(V_{L_{2}})_{1}$ has an
orthonormal basis:
$$x^{1}=\frac{1}{\sqrt{2}}\al(-1){\bf 1}, \ x^2= \frac{1}{\sqrt{2}}(e^{\al}+e^{-\al}), \
x^3=\frac{i}{\sqrt{2}}(e^{\al}-e^{-\al}).$$ Let $\tau_{i}\in
Aut(V_{L_{2}})$, $i=1,2,3$ be such that
$$
\tau_{1}(x^1,x^2,x^3)=(x^1,x^2,x^3)\left[\begin{array}{ccc}1& &\\
&-1&\\&&-1\end{array}\right],$$
$$
\tau_{2}(x^1,x^2,x^3)=(x^1,x^2,x^3)\left[\begin{array}{ccc}-1& &\\
&1&\\&&-1\end{array}\right],$$
$$
\tau_{3}(x^1,x^2,x^3)=(x^1,x^2,x^3)\left[\begin{array}{ccc}-1& &\\
&-1&\\&&1\end{array}\right].$$ Let $\sigma\in Aut(V_{L_{2}})$ be
such that
$$
\sigma(x^1,x^2,x^3)=(x^1,x^2,x^3)\left[\begin{array}{ccc}0&1 &0\\
0&0&-1\\-1&0&0\end{array}\right].$$ Then $\sigma$ and
$\tau_{i},i=1,2,3$ generate a finite subgroup of $Aut(V_{L_{2}})$
isomorphic to the alternating group $A_{4}$. We simply denote this
subgroup by $A_{4}$. It is easy to check that the subgroup $K$
generated by $\tau_{i}$, $i=1,2,3$ is a normal subgroup of $A_{4}$
of order $4$. Let
$$J= h(-1)^4{\bf 1} -2h(-3)h(-1){\bf 1} + \frac{3}{2}h(-2)^2{\bf
1}, \ \ E=e^{\be}+e^{-\be}$$ where $h=\frac{1}{\sqrt{2}}\al$,
$\be=2\al$. The following lemma comes from \cite{DG}.
\begin{lem}\label{l3.1} $V_{L_{2}}^K\cong V_{\Z\be}^{+}$, and $V_{\Z\be}^{+}$ is generated by $J$ and $E$.
Moreover, $(V_{L_{2}}^K)_4$ is four dimensional with a basis
$L(-2)^2\1, L(-4)\1, J, E.$
\end{lem}
By Lemma \ref{l3.1}, we have
$V_{L_{2}}^{A_{4}}=(V_{\Z\be}^+)^{\<\sigma\>}.$ A direct calculation
yields that
\begin{lem}\label{l3.2} We have
$$
\sigma(J)=-\frac{1}{2}J+\frac{9}{2}E, \ \
\sigma(E)=-\frac{1}{6}J-\frac{1}{2}E.
$$
\end{lem}

Let
\begin{equation}\label{e3.9}X^1=J-\sqrt{27}iE, \ X^2=J+\sqrt{27}iE.
\end{equation}
 Then it is easy to check that
\begin{equation}\label{e3.1}
\sigma(X^1)=\frac{-1+\sqrt{3}i}{2}X^1,
 \  \
\sigma(X^2)=\frac{-1-\sqrt{3}i}{2}X^2.
\end{equation}
It follows that $(V_{\Z\be}^+)^{\<\sigma\>}_4\subset L(1,0)$ and
$$(V_{\Z\be}^+)^{\<\sigma\>}=L(1,0)\bigoplus \sum_{n\geq 3}a_nL(1,n^2)$$
as a module for $L(1,0)$, where $a_n$ is the multiplicity of
$L(1,n^2)$ in $(V_{\Z\be}^+)^{\<\sigma\>}.$
 By (\ref{e3.1}) we immediately have for any $n\in\Z$,
$$
X^{1}_{n}X^2\in(V_{\Z\be}^+)^{\<\sigma\>}.
$$

For convenience, we call highest weight vectors for the Virasoro
algebra  primary vectors.  Note from \cite{DG} that $V_{\Z\be}^{+}$
contains two linearly independent primary vectors   $J$ and $E$ of
weight $4$ and one linearly independent primary vector of weight
$9.$ It is straightforward to compute that
$$
J_{3}J=-72L(-4){\bf 1}+336L(-2)^2{\bf 1}-60J, \
E_{3}E=-\frac{8}{3}L(-4){\bf 1}+\frac{ 112}{9}L(-2)^2{\bf
1}+\frac{20}{9}J$$ (cf. \cite{DJ3}). By Theorem \ref{2t1} and Lemma
\ref{l3.1}, we have for $n\in\Z$
$$
X^1_{n}X^2\in L(1,0)\oplus L(1,9)\oplus L(1,16).
$$
Note that $\sigma(E_{-2}J-J_{-2}E)=E_{-2}J-J_{-2}E.$ Since
$E_{-2}J-J_{-2}E\in M(1)\otimes e^{\beta}+M(1)\otimes e^{-\beta}$ we
see immediately that $E_{-2}J-J_{-2}E$ is a primary vector of weight
$9.$

The following lemma follows  from Theorem 3 in \cite{DM1} and
(\ref{e3.1}).
\begin{lem}\label{deco} We have decomposition
$$V_{\Z\be}^{+}=(V_{\Z\be}^{+})^{0}\oplus (V_{\Z\be}^{+})^{1}\oplus (V_{\Z\be}^{+})^{2},$$
where $(V_{\Z\be}^{+})^{0}=(V_{\Z\be}^{+})^{\<\sigma\>}$ is a simple
vertex operator algebra and $(V_{\Z\be}^{+})^{i}$ is the irreducible
$(V_{\Z\be}^{+})^{0}$-module generated by $X^i$ with lowest weight
$4$, $i=1,2$.
\end{lem}
 Set
\begin{eqnarray}
& u^0=-\dfrac{8}{3}L(-4){\bf 1}+\dfrac{112}{9}L(-2)^{2}{\bf
1}\\\label{a1} & u^1=-\dfrac{16}{9}L(-5){\bf
1}+\dfrac{112}{9}L(-3)L(-2){\bf 1}\\\label{a2} &
u^2=(-\dfrac{1856}{135}L(-6)-\dfrac{2384}{135}L(-4)L(-2)+\dfrac{1316}{135}L(-3)^2+\dfrac{1088}{135}L(-2)^3){\bf
1}\\\label{a3} &
u^3=(-\dfrac{464}{45}L(-7)-\dfrac{928}{45}L(-5)L(-2)+\dfrac{40}{9}L(-4)L(-3)+\dfrac{544}{45}L(-3)L(-2)^2){\bf
1}
\end{eqnarray}
\begin{eqnarray}
&v^2=(\dfrac{28}{75}L(-2)+\dfrac{23}{300}L(-1)^2)J,\\\label{a4}
&v^3=(\dfrac{14}{75}L(-3)+\dfrac{14}{75}L(-2)L(-1)-\dfrac{1}{300}L(-1)^3)J.\label{a5}
\end{eqnarray}
\begin{eqnarray}
&v^4=(\dfrac{28}{75}L(-2)+\dfrac{23}{300}L(-1)^2)E,\\\label{a6}
&v^5=(\dfrac{14}{75}L(-3)+\dfrac{14}{75}L(-2)L(-1)-\dfrac{1}{300}L(-1)^3)E.\label{a7}
\end{eqnarray}
By Lemma 2.5 of \cite{DJ3}, we have the following lemma.
\begin{lem}\label{l3.6} Let $E$ and $J$ be as before. Then
\begin{eqnarray*}
& E_{3}E =u^0+\dfrac{20}{9}J, \ J_{3}J=27u^0-60J, \ J_{3}E=60E,\\
& E_{2}E=u^1+\dfrac{10}{9}L(-1)J, \ J_{2}J=27u^1-30L(-1)J, \ J_{2}E=30L(-1)E,\\
& E_{1}E=u^2+\dfrac{20}{9}v^2, \ J_{1}J=27u^2-60v^2, \ J_{1}E=60v^4,\\
& E_{0}E=u^3+\dfrac{20}{9}v^3, \ J_{0}J=27u^3-60v^3, \ J_{0}E=60v^5.
\end{eqnarray*}
\end{lem}
Using Lemma \ref{l3.6}, one can check directly that
\begin{equation}\label{action on J}
(J_{-2}E-E_{-2}J)_{8}J=-10800E, \ (J_{-2}E-E_{-2}J)_{8}E=400J.
\end{equation}
As a result, we have
\begin{equation}\label{e3.2}
(J_{-2}E-E_{-2}J)_{8}X^{1}=-400\sqrt{27}iX^1, \
(J_{-2}E-E_{-2}J)_{8}X^2=400\sqrt{27}iX^2.
\end{equation}
By (\ref{e3.2}), we immediately know that $J_{-2}E-E_{-2}J$ is a
non-zero primary vector of weight $9.$ Recall from \cite{DG} that
$V_{\Z\be}^{+}$ has one primary vector of weight $9$ up to a
constant. A direct calculation yields that
\begin{lem}\label{l3.4} The vector
\begin{eqnarray*}
u^{(9)}=&
-\dfrac{1}{\sqrt{2}}(15h(-4)h(-1)+10h(-3)h(-2)+10h(-2)h(-1)^3)\otimes
(e^{\be}+e^{-\be})\\
& + (6h(-5)+10h(-3)h(-1)^2+\dfrac{15}{2}h(-2)^2h(-1)+h(-1)^5)\otimes
(e^{\be}-e^{-\be})
\end{eqnarray*}
is a non-zero primary vector of weight $9$ and
$u^{(9)}\in\C(J_{-2}E-E_{-2}J)$.
\end{lem}
Note from \cite{L1} that there is a non-degenerate symmetric
invariant bilinear form $(\cdot,\cdot)$ on $V_{\Z\beta}^+.$ The next
lemma gives a relation between  $u^{(9)}$ and $J_{-2}E-E_{-2}J.$
\begin{lem}\label{relation} We have
$$
J_{-2}E-E_{-2}J=-2\sqrt{2}u^{(9)},
$$
$$
(u^{(9)},u^{(9)})=5400.
$$
\end{lem}
\pf By Lemma \ref{l3.4} and (\ref{action on J}), we have
$u^{(9)}_{8}E\in\C J.$
 Denote $F=e^{\be}-e^{-\be}$. Note that
$$ (V_{\Z\be}^{+})_{4}=\C h(-3)h(-1){\bf 1}\oplus \C h(-2)^{2}{\bf
1}\oplus \C h(-1)^4{\bf 1}\oplus \C E.
$$
%and $$J=h(-1)^4{\bf 1}-2h(-3)h(-1){\bf 1} + \frac{3}{2}h(-2)^2{\bf 1}.$$
Let $W_{4}$ be the subspace of $(V_{\Z\be}^{+})_{4}$ linearly
spanned by $E$, $h(-3)h(-1){\bf 1}$ and $h(-2)^2{\bf 1}$. Then
$$
h(-1)^4{\bf 1}\equiv J\ (\mod\ W_{4}).$$ Furthermore,  we have
\begin{align*}
(h(-4)h(-1)E)_{8}E \equiv &
\sum\limits_{i=0}^{\infty}(-1)^{i+1}\left(\begin{array}{c}-4\\i\end{array}\right)(h(-1)E)_{4-i}h(i)E
(\mod\ W_{4})\\
\equiv & -\sqrt{8}(h(-1)E)_{4}F (\mod\ W_{4})\\
\equiv & -\sqrt{8}
\sum\limits_{i=0}^{\infty}(-1)^{i}\left(\begin{array}{c}-1\\i\end{array}\right)(h(-1-i)E_{4+i}+E_{3-i}h(i))F (\mod\ W_{4})\\
\equiv & -\sqrt{8}h(-1)E_{4}F-8E_{3} E\,(\mod \ W_{4}).
\end{align*}
 Similarly, \begin{align*}
 (h(-3)h(-2)E)_{8}E\equiv & -8E_{3}E\,
(\mod\ W_{4}),
\\
(h(-2)h(-1)^3E)_{8}E\equiv &
-\sqrt{8}h(-1)^3E_{6}F-24h(-1)^2E_{5}E\\
&-24\sqrt{8}h(-1)E_{4}F-64E_{3}E\, (\mod\ W_{4}), \\
(h(-5)F)_{8}E\equiv & \sqrt{8}F_{3}F\ (\mod\ W_{4}),
\\
(h(-3)h(-1)^{2}F)_{8}E\equiv &
\sqrt{8}h(-1)^2F_{5}F+16h(-1)F_{4}E+8\sqrt{8}F_{3}F\ (\mod\
W_{4}).\,
\\
(h(-2)^2h(-1)F)_{8}E\equiv & 8h(-1)F_{4}E+8\sqrt{8}F_{3}F\ (\mod\
W_{4}),
\\
(h(-1)^5F)_{8}E\equiv &
5\sqrt{8}h(-1)^{4}F_{7}F+80h(-1)^3F_{6}E+80\sqrt{8}h(-1)^2F_{5}F\\
&+320h(-1)F_{4}E+64\sqrt{8}F_{3}F\ (\mod\ W_{4}).
\end{align*}
It is then easy to check that
$$u^{(9)}_{8}E=-100\sqrt{2}h(-1)^4{\bf 1}\ (\mod\ W_{4}).$$
This implies that
$$
u^{(9)}_{8}E=-100\sqrt{2}J.
$$
Then by (\ref{action on J}),
$$
J_{-2}E-E_{-2}J=-2\sqrt{2}u^{(9)}.
$$
Note that
$$
(J_{-2}E-E_{-2}J,J_{-2}E)=(J_{8}(J_{-2}E-E_{-2}J),E)=-((J_{-2}E-E_{-2}J)_{8}J,E)
$$
and
$$
(J_{-2}E-E_{-2}J,E_{-2}J)=(E_{8}(J_{-2}E-E_{-2}J),J)=(-(J_{-2}E-E_{-2}J)_{8}E,J).
$$
Since
$$(E,E)=2,\ (J,J)=54,$$
(see \cite{DJ2}) it follows from (\ref{action on J}) that
\begin{equation}\label{relation1 of two vectors}
(J_{-2}E-E_{-2}J,J_{-2}E-E_{-2}J)=43200
\end{equation}
and
$$(u^{(9)},u^{(9)})=5400.$$
The proof is complete. \qed

\section{$C_{2}$-cofiniteness and rationality of $V_{L_{2}}^{A_{4}}$
}
\def\theequation{4.\arabic{equation}}
\setcounter{equation}{0} The $C_{2}$-cofiniteness and rationality of
$V_{L_{2}}^{A_{4}}$ is established in this section. The proof
involves some very hard computations.

By Lemma \ref{l3.2}, we have
$$
J_{-9}J+27E_{-9}E\in (V_{\Z\be}^{+})^{\<\sigma\>}.
$$
Then it is clear that
\begin{equation}\label{eq4.1}
J_{-9}J+27E_{-9}E=x^0+X^{(16)}+27(e^{2\be}+e^{-2\be}),
\end{equation}
where $x^0\in L(1,0)$, and  $X^{(16)}$  is a non-zero primary
element of weight 16 in $M(1)^{+}$. Denote
\begin{equation}\label{weight16}
u^{(16)}=X^{(16)}+27(e^{2\be}+e^{-2\be}).
\end{equation}
Then $u^{(16)}\in(V_{\Z\be}^{+})^{\<\sigma\>}$ is a non-zero primary
vector of weight 16.
\begin{lem}\label{lweight16}We have the following:
$$
u^{(9)}_{1}u^{(9)}-58800u^{(16)}\in L(1,0).$$
\end{lem}
 \pf  Denote $E^2=e^{2\be}+e^{-2\be}$. By Theorem \ref{2t1} and the skew-symmetry, we may assume that
$$
u^{(9)}_{1}u^{(9)}=v+cu^{(16)},
$$
for some $v\in L(1,0)$ and $c\in\C$. To determine $c$ we just need
to consider $(u^{(9)}_{1}u^{(9)},E^2)$ by (\ref{weight16}). Recall
that
\begin{eqnarray*}
u^{(9)}=&
-\dfrac{1}{\sqrt{2}}(15h(-4)h(-1)+10h(-3)h(-2)+10h(-2)h(-1)^3)\otimes E\\
& + (6h(-5)+10h(-3)h(-1)^2+\dfrac{15}{2}h(-2)^2h(-1)+h(-1)^5)\otimes
F,
\end{eqnarray*}
where $F=e^{\be}-e^{-\be}$. To calculate $((h(-4)h(-1)\otimes
E)_{1}(h(-4)h(-1)\otimes E), E^2)$, we only need to consider the
coefficient of the monomial $E^2$ in $(h(-4)h(-1)\otimes
E)_{1}(h(-4)h(-1)\otimes E)$. Then direct calculation yields that
$$
 ((h(-4)h(-1)\otimes E)_{1}(h(-4)h(-1)\otimes )E),
 E^2)=(972E^2,E^2).$$
Calculations for other monomials are similar. For example,
$$
((h(-3)h(-2)\otimes E)_{1}(h(-2)h(-1)^3\otimes )E),
 E^2)=(304E^2,E^2).
$$
Then one can check that
$$
(u^{(9)}_{1}u^{(9)},E^2)=(1587600E^2,E^2).$$ It follows that
$c=58800$.
 \qed

\begin{lem}\label{generators} The following hold:
(1) \ $(V_{\Z\be}^{+})^{\<\sigma\>}$ is generated by $u^{(9)}$.

(2) \  $(V_{\Z\be}^{+})^{<\sigma>}$ is linearly spanned by
$$
L(-m_{s})\cdots L(-m_{1})u^{(9)}_{n}u^{(9)}, \ L(-m_{s})\cdots
L(-m_{1})w^{p}_{-k_{p}}\cdots w^{1}_{-k_{1}}w,$$ where
$w,w^1,\cdots,w^p\in\{u^{(9)},u^{(16)}\}$, $k_{p}\geq \cdots\geq
k_{1}\geq 2$, $n\in\Z$, $m_{s}\geq\cdots \geq m_{1}\geq 1$, $s,p\geq
0$.

\end{lem}

\pf  By Lemma \ref{relation}, $\omega$ can be generated by
$u^{(9)}$. It follows from \cite{DGR} that
$(V_{\Z\be}^{+})^{\<\sigma\>}$ is generated by $u^{(9)}$ and
$u^{(16)}$. Then (1)  follows from Lemma \ref{lweight16}.

 By (3.2) in \cite{A1} and (3.3) in \cite{A3}, we have
\begin{equation}\label{Hei1}
M(1,2\sqrt{2}m)=\bigoplus_{p=0}^{\infty}L(1,(2m+p)^2),
\end{equation}
\begin{equation}\label{Hei2}
V_{\Z\be}^{+}=M(1)^{+}\bigoplus(\bigoplus_{m=1}^{\infty}M(1,2\sqrt{2}m))=M(1)^{+}\bigoplus(\bigoplus_{m=1}^{\infty}(\bigoplus_{p=0}^{\infty}L(1,(2m+p)^2)).
\end{equation}
 By (\ref{Hei2}) the subspace $U^1$ linearly spanned by  primary
elements of weight 16 in $V_{\Z\be}^{+}$ is three dimensional.
Obviously $U^1$ is invariant under $\sigma$.  Note that
$e^{2\be}+e^{-2\be}\in U^1$. Consider the $M(1)^{+}$-submodule $W$
of $V_{\Z\be}^{+}$ generated by $e^{2\be}+e^{-2\be}$. If
$e^{2\be}+e^{-2\be}\in (V_{\Z\be}^{+})^{<\sigma>}$, then by the
fusion rule of $M(1)^{+}$ (also see \cite{DN2}), $J\in W\cdot
W=\<u_{n}v|u,v\in W, n\in\Z\>.$ So $J\in
(V_{\Z\be}^{+})^{<\sigma>},$ which contradicts with Lemma
\ref{l3.2}. This implies that $\sigma$ has an eigenvector in $U^1$
with eigenvalue not equal to 1. Since $\sigma^3=1$ and $U^1$ is a
vector space over  $\C$, it follows that both
$\frac{-1+\sqrt{3}i}{2}$ and $\frac{-1-\sqrt{3}i}{2}$ occur as
eigenvalues of $\sigma$ on $U_1.$
 Recall that $u^{(16)}\in
(V_{\Z\be}^{+})^{<\sigma>}$ is a non-zero primary element of weight
16. So we immediately have
$$
(V_{\Z\be}^{+})^{<\sigma>}=L(1,0)\bigoplus L(1,9)\bigoplus
L(1,16)\bigoplus(\sum\limits_{n\geq 5}a_{n}L(1,n^2)).
$$

Let $U^2$ be the subspace of primary vectors of weight 25 in
$V_{\Z\be}^{+}$. By (\ref{Hei2}), $\dim U^2=2$. Consider the
$M(1)^{+}$-submodule $W$ of $V_{\Z\be}^{+}$ generated by
$e^{2\be}+e^{-2\be}$ again. By (\ref{Hei1}), there is a non-zero
primary element $w=x\otimes (e^{2\be}+e^{-2\be})+y\otimes
(e^{2\be}-e^{-2\be})$ of weight 25 in $W$ for some $x\in M(1)^{+}$
and $y\in M(1)^{-}$. Let $U^{(9)}$ and $U^{(16)}$ be the
$L(1,0)$-submodules of $(V_{\Z\be}^{+})^{<\sigma>}$ generated by
$u^{(9)}$ and $u^{(16)},$ respectively. Then by Part (1) and the
skew-symmetry, any element of weight 25 in
$(V_{\Z\be}^{+})^{<\sigma>}$ is a linear combination of elements in
$L(1,0)\oplus U^{(9)}\oplus U^{(16)}$ and $U^{(9)}\cdot
U^{(16)}=<u_{n}v|u\in U^{(9)}, v\in U^{(16)}>$. By Lemma \ref{l3.4}
and  (\ref{weight16}), elements in $U^{(9)}\cdot U^{(16)}$ have the
forms: $u\otimes (e^{\be}+e^{-\be})+v\otimes (e^{\be}-e^{-\be})$,
where $u\in M(1)^+$ and $v\in M(1)^-$. So we know that $w\notin
(V_{\Z\be}^{+})^{<\sigma>}$. This proves that $\sigma|_{U^2}$ has
eigenvalues not equal to 1. Since $\sigma^3=1$ and $\dim_{\C}U^2=2$,
it follows that $\sigma|_{U^2}$ has two eigenvalues
$\frac{-1+\sqrt{3}i}{2}$ and $\frac{-1-\sqrt{3}i}{2}$. So we
immediately have
\begin{equation}\label{adde}
(V_{\Z\be}^{+})^{<\sigma>}=L(1,0)\bigoplus L(1,9)\bigoplus
L(1,16)\bigoplus(\sum\limits_{n\geq 6}a_{n}L(1,n^2)).
\end{equation}
A proof similar to that of Lemma 4.3 in \cite{DN2} gives (2) with
the help of (1), (\ref{adde}) and Theorem \ref{2t1}.
 \qed

\vskip 0.3cm \begin{lem}\label{weight20} We have
\begin{align*}
u^{(9)}_{-3}u^{(9)}=& s^1+\frac{162770}{99}L(-4)u^{(16)}+
\frac{5204015}{1584}L(-3)L(-1)u^{(16)}+\frac{14760}{11}L(-2)^{2}u^{(16)}\\&
+\frac{1154225}{792}L(-2)L(-1)^2u^{(16)}+\frac{354895}{3168}L(-1)^{4}u^{(16)},
\end{align*}
\begin{align*}
u^{(9)}_{-5}u^{(9)}=&
s^2-\frac{653871670}{63063\cdot 27}L(-6)u^{(16)}+\frac{3303230375}{2018016\cdot 27}L(-5)L(-1)u^{(16)}\\
&
+\frac{489993820}{63063\cdot 27}L(-4)L(-2)u^{(16)}+\frac{69658220}{9009\cdot 27}L(-3)^2u^{(16)}\\
& +\frac{346772585}{42042\cdot
27}L(-4)L(-1)^2u^{(16)}+\frac{3338006885}{168168\cdot
27}L(-3)L(-2)L(-1)u^{(16)}\\
&+ \frac{19408720}{7007\cdot 27}L(-2)^3u^{(16)}
+\frac{14067649205}{4036032\cdot 27}L(-3)L(-1)^3u^{(16)}\\
&+\frac{1055175305}{252252\cdot
27}L(-2)^2L(-1)^2u^{(16)}+\frac{1185150565}{2018016\cdot
27}L(-2)L(-1)^4u^{(16)}\\
&+\frac{119070745}{8072064\cdot 27}L(-1)^6u^{(16)},
\end{align*}
where $s^1,s^2\in L(1,0)$.
\end{lem}
\pf  By  Theorem \ref{2t1} and the skew-symmetry, we may assume that
$$
u^{(9)}_{-3}u^{(9)}=s^1+y^1, \ u^{(9)}_{-5}u^{(9)}=s^2+y^2,
$$
where $s^1,s^2\in L(1,0)$, $y^1,y^2\in U^{(16)}\cong L(1,16)$ which
is an $L(1,0)$-submodule of $(V_{\Z\be}^{+})^{\<\sigma\>}$ generated
by $u^{(16)}$. Then we may assume that
\begin{align*}
y^1=&
a_{1}L(-4)u^{(16)}+a_{2}L(-3)L(-1)u^{(16)}+a_{3}L(-2)^2u^{(16)}\\
& +a_{4}L(-2)L(-1)^2u^{(16)}+a_{5}L(-1)^4u^{(16)}\\
= & \sum\limits_{i=1}^{5}a_{i}w^{i}.
\end{align*}
To determine $a_{i},1\leq i\leq 5$, we consider
$(u^{(9)}_{-3}u^{(9)},w^i)$, $(w^i,w^j)$, $i,j=1,2,\cdots,5$. Then
by Lemma \ref{lweight16} and direct calculation, we have
$$
\left[\begin{array}{ccccc}133& 224& 387& 576& 1920\\ 224& 3328& 480&
10560& 49920\\ 387& 480&
   17673/2& 13152& 57600\\ 576& 10560& 13152& 162336&
   1267200\\ 1920& 49920& 57600& 1267200& 30159360\end{array}\right]\left[\begin{array}{c} a_{1}\\
 a_{2}\\a_{3}\\a_{4}\\a_{5}\end{array}\right]=58800\left[\begin{array}{c} 43\\
 560\\675\\7344\\93024\end{array}\right].
$$
We get that
$$a_{1}=\frac{162770}{99}, \ a_{2}=\frac{5204015}{1584},$$$$  a_{3}= \frac{14760}{11},\  a_{4}=\frac{1154225}{792}, \
a_{5}=\frac{354895}{3168}.$$ The  first formula follows. The proof
for the second  one is similar. We omit it.
 \qed

 \vskip 0.3cm

  Let $v$ be any element in $V_{\Z\be}^{+}$ of weight $m\leq 22$. Then
$v$ is a linear combination of an element in $V^{(4)}\oplus
V^{(16)}$ and elements in $M(1)^{+}$ having the form
$h(-n_{t})\cdots h(-n_{1}){\bf 1}$ such that $n_{t}\geq \cdots\geq
n_{1}\geq 1$ and $\sum\limits_{i=1}^t n_{i}=m$, where $V^{(4)}$ and
$V^{(16)}$ are $M(1)^{+}$-submodules of $V_{\Z\be}^{+}$ generated by
$E$ and $E^2$ respectively. We denote by $c(v)$ the coefficient of
the monomial $h(-1)^{m}{\bf 1}$ in the linear combination.  Then
 we have the following lemma.
 \begin{lem}\label{coefficient}
 $$c(u^{(9)}_{-3}u^{(9)})=-\frac{447232}{19\cdot 17\cdot 11\cdot 7^2\cdot 5^2\cdot 3},$$
 $$
 c(u^{(9)}_{-5}u^{(9)})=-\frac{328099328}{19\cdot 17\cdot
 13\cdot 11^2\cdot 7^3\cdot 5^2\cdot 3^6}.
 $$
 \end{lem}
\pf  Let $k\in2\Z+1$. We consider $c(u^{(9)}_{-k}u^{(9)})$. By a
direct but long  calculation, we have
\begin{align*}
c(u^{(9)}_{-k}u^{(9)})=& -2700 c(E_{-k-10}E) -
13500c(h(-1)^2E_{-k-8}E) - 18000\sqrt{2}c(h(-1)^3E_{-k-7}F)\\
& -31500c(h(-1)^4E_{-k-6}E)-15300\sqrt{2}c(h(-1)^5E_{-k-5}F)\\
& -9060c(h(-1)^6E_{-k-4}E)-1620\sqrt{2}c(h(-1)^7E_{-k-3}F) \\
& -345c(h(-1)^8E_{-k-2}E)-20\sqrt{2}c(h(-1)^9E_{-k-1}F)
-c(h(-1)^{10}E_{-k}E)
\end{align*}
Note that for $m\in2\Z+1$, $n\in2\Z$, $m,n\leq 7$,
$$
c(E_{m}E)=\frac{2\cdot (\sqrt{8})^{7-m}}{(7-m)!}, \
c(E_{n}F)=-\frac{2\cdot (\sqrt{8})^{7-n}}{(7-n)!}.$$ Let $k=-3,
k=-5$ respectively, we then get the lemma. \qed

\vskip 0.3cm As defined in \cite{Z}, a vertex operator algebra $V$
is called $C_{2}$-cofinite, if $V/C_{2}(V)$ is finite-dimensional,
where $C_{2}(V)=span_{\C}\{u_{-2}v|u,v\in V\}$. The following lemma
comes from \cite{Z}.

\begin{lem}\label{prop of cofi}
(1) $L(-1)u\in C_{2}(V)$ for $u\in V$;

(2) $u_{-k}v\in C_{2}(V)$, for $u,v\in V$ and $k\geq 2$;

(3) $u_{-1}v\in C_{2}(V)$, for $u\in V, v\in C_{2}(V)$.
\end{lem}
We are now in a position to state the main result of this section.
\begin{theorem}\label{rationality}
 $(V_{\Z\be}^{+})^{\<\sigma\>}$ is $C_{2}$-cofinite and rational.
\end{theorem}
\pf Let $s^1,s^2\in L(1,0)$ be the same as in Lemma \ref{weight20}.
Then $s^1$ and $s^2$ are linear combinations of linearly independent
elements having the forms $L(-m_{s})\cdots L(-m_{1}){\bf 1}$ and
$L(-n_{t})\cdots L(-n_{1}){\bf 1}$ respectively such that
$m_{s}\geq\cdots\geq  m_{1}\geq 2$, $n_{t}\geq \cdots\geq n_{1}\geq
2$ and $\sum\limits_{i=1}^{s}m_{i}=20$,
$\sum\limits_{i=1}^{t}n_{i}=22$. Assume the coefficients of
$L(-2)^{10}{\bf 1}$ and $L(-2)^{11}{\bf 1}$ in the two linear
combinations are $a_{1}$ and $a_{2}$ respectively. Then by Lemma
\ref{prop of cofi},
$$
s^1-a_{1}L(-2)^{10}{\bf 1},\  s^2-a_{2}L(-2)^{11}{\bf 1}\in
C_{2}((V_{\Z\be}^{+})^{\<\sigma\>}).
$$
Further, by Lemma \ref{prop of cofi} and Lemma \ref{weight20}, we
have
$$
s^1+\frac{14760}{11}L(-2)^{2}u^{(16)}, \
s^2+\frac{19408720}{7007\cdot 27}L(-2)^3u^{(16)}\in
C_{2}((V_{\Z\be}^{+})^{\<\sigma\>}).$$ So
$$
a_{1}L(-2)^{10}{\bf 1}+\frac{14760}{11}L(-2)^{2}u^{(16)}, \
a_{2}L(-2)^{11}{\bf 1}+\frac{19408720}{7007\cdot
27}L(-2)^3u^{(16)}\in C_{2}((V_{\Z\be}^{+})^{\<\sigma\>}).$$ Thus by
Lemma \ref{prop of cofi} \begin{equation}\label{cofi3}
a_{1}L(-2)^{11}{\bf 1}+\frac{14760}{11}L(-2)^{3}u^{(16)}, \
a_{2}L(-2)^{11}{\bf 1}+\frac{19408720}{7007\cdot
27}L(-2)^3u^{(16)}\in C_{2}((V_{\Z\be}^{+})^{\<\sigma\>}).
\end{equation}

On the other hand, note from the definition of $L(-2){\bf 1}$ that
$c(L(-2)^k\1)=2^k.$ This implies that
$$
c(u^{(9)}_{-3}u^{(9)})=\frac{1}{2^{10}}a_{1}+\frac{1}{4}\cdot\frac{14760}{11}c(X^{(16)}),
$$
$$
c(u^{(9)}_{-5}u^{(9)})
=\frac{1}{2^{11}}a_{2}+\frac{1}{8}\cdot\frac{19408720}{7007\cdot
27}c(X^{(16)}).$$ So by Lemma \ref{coefficient},
\begin{equation}\label{cofi1}
\frac{1}{2^{10}}a_{1}+\frac{1}{4}\frac{14760}{11}c(X^{(16)})=-\frac{447232}{19\cdot
17\cdot 11\cdot 7^2\cdot 5^2\cdot 3},
\end{equation}
\begin{equation}\label{cofi2}
\frac{1}{2^{11}}a_{2}+\frac{1}{8}\frac{19408720}{7007\cdot
27}c(X^{(16)})=-\frac{328099328}{19\cdot 17\cdot
 13\cdot 11^2\cdot 7^3\cdot 5^2\cdot 3^6}.
\end{equation}
 If $$a_{1}/a_{2}=\frac{14760}{11}/\frac{19408720}{7007\cdot
27},$$ then by (\ref{cofi1}) and (\ref{cofi2}), we have
$$
\frac{-\frac{447232}{2\cdot19\cdot 17\cdot 11\cdot 7^2\cdot 5^2\cdot
3}}{-\frac{328099328}{19\cdot 17\cdot
 13\cdot 11^2\cdot 7^3\cdot 5^2\cdot 3^6}}=\frac{\frac{14760}{11}}{\frac{19408720}{7007\cdot
27}}.
$$
But
$$
\frac{-\frac{447232}{2\cdot 19\cdot 17\cdot 11\cdot 7^2\cdot
5^2\cdot 3}}{-\frac{328099328}{19\cdot 17\cdot
 13\cdot 11^2\cdot 7^3\cdot 5^2\cdot 3^6}}=\frac{32688117}{2563276}\neq \frac{6346431}{485218}=\frac{\frac{14760}{11}}{\frac{19408720}{7007\cdot 27}}.
$$
This means that
$$
a_{1}/a_{2}\neq \frac{14760}{11}/\frac{19408720}{7007\cdot 27}.$$ By
(\ref{cofi3}), we have
$$
L(-2)^{11}{\bf 1}, L(-2)^{3}u^{(16)}\in
C_{2}((V_{\Z\be}^{+})^{\<\sigma\>}).
$$
Then it follows from Lemma \ref{generators} that
$(V_{\Z\be}^{+})^{\<\sigma\>}$ is $C_{2}$-cofinite. Since
$V_{\Z\be}^{+}$ is rational and $(V_{\Z\be}^{+})^{\<\sigma\>}$ is
self-dual, it follows that $(V_{\Z\be}^{+})^{\<\sigma\>}$ satisfies
the Hypothesis I in \cite{M2}. Then by Corollary 7 in \cite{M2},
$(V_{\Z\be}^{+})^{\<\sigma\>}$ is rational. \qed

\section{Classification and construction of irreducible modules of $(V_{\Z\be}^{+})^{<\sigma>}$}
\def\theequation{5.\arabic{equation}}
\setcounter{equation}{0}

 We  will first  construct all the irreducible $\sigma^i$-twisted modules of $V_{\Z\be}^{+}$, $i=1,2$. We  have the following
lemma.
\begin{lem}\label{twisted-modules}
There are at most two inequivalent irreducible $\sigma$-twisted
modules of $V_{\Z\be}^{+}$.
\end{lem}

\pf Let $(W,Y)$ be an irreducible $V_{\Z\be}^{+}$-module. Define a
linear map
$$Y^{\sigma}: V_{\Z\be}^{+}\rightarrow (\End
W)[[z,z^{-1}]]$$
 by
$$Y^{\sigma}(u,z)w=Y(\sigma^{-1}(u),z)w$$
where $u\in V_{\Z\be}^{+}$, $w\in W$. Recall from \cite{DLM1} that
$(W,Y^{\sigma})$ is still an irreducible module of $V_{\Z\be}^{+}$,
which we denote by $W^{\sigma}$. As in \cite{DLM1}, if $W\cong
W^{\sigma}$, we say $W$ is stable under $\sigma$. Recall from
\cite{DN2} that all the irreducible modules of $V_{\Z\be}^{+}$ are
$$
V_{\Z\be}^{\pm}, \  V_{\Z\be+\frac{r}{8}\be} \ (1\leq r\leq 3), \
V_{\Z\be+\frac{\be}{2}}^{\pm}, \ V_{\Z\be}^{T_{1},\pm}, \
V_{\Z\be}^{T_{2},\pm}
$$
with the following tables
 \vskip 5ex
\begin{center}
\begin{tabular}{|c|c|c|c|c|c|c|c|}
\hline
&$V_{\Z\be}^{+}$&$V_{\Z\be}^{-}$&$V_{\Z\be+\frac{1}{8}\be}$&$V_{\Z\be+\frac{1}{4}\be}$&$V_{\Z\be+\frac{3}{8}\be}$&
$V_{\Z\be+\frac{\be}{2}}^+$&$V_{\Z\be+\frac{\be}{2}}^-$\\
\hline $\omega$&$0$&$1$&$\frac{1}{16}$&$\frac{1}{4}$&$\frac{9}{16}$&$1$&$1$\\
\hline
$E$&$0$&$0$&$0$&$0$&$0$&$1$&$-1$\\
\hline $J$&$0$&$-6$&$-\frac{3}{64}$&$0$&$\frac{45}{64}$&$3$&$3$\\
\hline
\end{tabular}
\end{center}
\vskip 5ex
\begin{center}
\begin{tabular}{|c|c|c|c|c|}
\hline
&$V_{\Z\be}^{T_{1},+}$&$V_{\Z\be}^{T_{1},-}$&$V_{\Z\be}^{T_{2},+}$&$V_{\Z\be}^{T_{2},-}$\\
\hline $\omega$&$1/16$&$9/16$&$1/16$&$9/16$\\
\hline
$E$&$1/128$&$-15/128$&$-1/128$&$15/128$\\
\hline $J$&$3/128$&$-45/128$&$3/128$&$-45/128$\\
\hline
\end{tabular}
\end{center}
\vskip 5ex

It is easy to check that
$$
V_{\Z\be}^{+}\cong (V_{\Z\be}^{+})^{\sigma},  \
(V_{\Z\be+\frac{\be}{4}})^{\sigma}\cong V_{\Z\be+\frac{\be}{4}},
$$$$ (V_{\Z\be}^{-})^{\sigma}\cong V_{\Z\be+\frac{\be}{2}}^-,  \
(V_{\Z\be+\frac{\be}{2}}^-)^{\sigma}\cong V_{\Z\be+\frac{\be}{2}}^+,
\ (V_{\Z\be+\frac{\be}{2}}^+)^{\sigma}\cong V_{\Z\be}^{-},
$$
$$
(V_{\Z\be+\frac{\be}{8}})^{\sigma}\cong V_{\Z\be}^{T_{2},+}, \
(V_{\Z\be}^{T_{2},+})^{\sigma}\cong V_{\Z\be}^{T_{1},+}, \
(V_{\Z\be}^{T_{1},+})^{\sigma}\cong V_{\Z\be+\frac{\be}{8}},
$$
$$
(V_{\Z\be+\frac{3\be}{8}})^{\sigma}\cong V_{\Z\be}^{T_{2},-}, \
(V_{\Z\be}^{T_{2},-})^{\sigma}\cong V_{\Z\be}^{T_{1},-}, \
(V_{\Z\be}^{T_{1},-})^{\sigma}\cong V_{\Z\be+\frac{3\be}{8}}.
$$
Then the lemma follows from \cite {A2}, \cite{Y} and Theorem 10.2 in
\cite{DLM1}. \qed

Next we will prove that there are exactly  two inequivalent
irreducible $\sigma$-twisted $V_{\Z\be}^{+}$-modules. We first
construct irreducible $\sigma$-twisted  $V_{L_{2}}$-modules. Let
$x^i,i=1,2,3$ be defined as in Section 3. Set
$$
h'=\dfrac{1}{3\sqrt{6}}(x^1+x^2-x^3), $$$$
y^1=\frac{1}{\sqrt{3}}(x^1+\dfrac{-1+\sqrt{3}i}{2}x^2+\dfrac{1+\sqrt{3}i}{2}x^3),
$$$$ y^2=\frac{1}{\sqrt{3}}(x^1+\dfrac{-1-\sqrt{3}i}{2}x^2+\dfrac{1-\sqrt{3}i}{2}x^3).
$$
Then
$$
L(n)h'=\delta_{n,0}h', \ h'(n)h'=\dfrac{1}{18}\delta_{n,1}\1, \
n\in\Z,
$$
$$
 h'(0)y^1=\frac{1}{3}y^1, \ h'(0)y^2=-\frac{1}{3}y^2, \
 y^1(0)y^2=6h'.
$$
 It follows that $h'(0)$ acts semisimply  on $V_{L_{2}}$ with
 rational eigenvalues. So $e^{2\pi ih'(0)}$ is an automorphism of
 $V_{L_{2}}$ (see \cite {Li2}, \cite{DG}, etc.). Since
$$
e^{2\pi ih'(0)}h'=h', \ e^{2\pi
ih'(0)}y^1=\dfrac{-1+\sqrt{3}i}{2}y^1, \ e^{2\pi
ih'(0)}y^2=\dfrac{-1-\sqrt{3}i}{2}y^2,
$$
it is easy to see that $$e^{2\pi ih'(0)}=\sigma.$$ Let
$$
\Delta(h',z)=z^{h'(0)}\exp(\sum\limits_{k=1}^{\infty}\dfrac{h'(k)}{-k}(-z)^{-k}),$$
and
$$ W^1=V_{L_{2}}, \ W^{2}=V_{L_{2}+\frac{\al}{2}}.$$
Then $W^1$ and $W^2$ are all the irreducible $V_{L_{2}}$-modules and
 $$
 W^1(0)=\C\1, \ W^2(0)=\C e^{\frac{\al}{2}}\bigoplus \C
 e^{-\frac{\al}{2}}.
$$
Let
$$
w^1=e^{\frac{\al}{2}}+\frac{(\sqrt{3}-1)(1+i)}{2}e^{-\frac{\al}{2}},$$
$$
w^2=\frac{1}{\sqrt{2}}[(\sqrt{3}-1)e^{\frac{\al}{2}}-(1+i)e^{-\frac{\al}{2}}].
$$
Then $W^2=\C w^1\bigoplus \C w^2$ and
$$
h'(0)w^1=\frac{1}{6}w^1, \ h'(0)w^2=-\frac{1}{6}w^2,
$$
$$
y^1(0)w^1=0, \ y^1(0)w^2=w^1, y^2(0) w^1=w^2.
$$

From \cite{Li2}, we have the following lemma.
\begin{lem}\label{twisted1}
$(W^{i,T_{1}},Y_{\sigma}(\cdot,z))=(W^i, Y(\Delta(h',z)\cdot,z))$
are irreducible $\sigma$-twisted modules of $V_{L_{2}}, i=1,2.$
\end{lem}
 Direct calculation yields that
\begin{equation}\label{virasoro element1}
\Delta(h',z)L(-2)\1=L(-2)\1+z^{-1}h'(-1)\1+\dfrac{1}{36}z^{-2}\1,
\end{equation}
\begin{equation}\label{s-element1}
Y_{\sigma}(h',z)=Y(h'+\frac{1}{18}z^{-1},z),
\end{equation}
\begin{equation}\label{u-nilpotent1}
Y_{\sigma}(y^1,z)=z^{\frac{1}{3}}Y(y^1,z),
\end{equation}
\begin{equation}\label{l-nilpotent1}
Y_{\sigma}(y^2,z)=z^{-\frac{1}{3}}Y(y^2,z).
\end{equation}
To distinguish the components of $Y(u,z)$ from those of
$Y_{\sigma}(u,z)$ we consider the following expansions
$$Y_{\sigma}(u,z)=\sum\limits_{n\in\Z+\frac{r}{3}}u_{n}z^{-n-1}, \
 Y(u,z)=\sum\limits_{n\in\Z}u(n)z^{-n-1},
 $$
where $u\in V_{L_{2}}$ such that $\sigma(u)=e^{-\frac{2r\pi
i}{3}}u.$
  By
 (\ref{s-element1})-(\ref{l-nilpotent1}) and direct calculation, we have the following
 lemma.
 \begin{lem}\label{twisted-module1} Write $W^{i,T_{1}}=\oplus_{n\in\frac{1}{3}\Z_+}W^{i,T_{1}}(n)$ as admissible
 $\sigma$-twisted module. Then
$$ W^{1,T_{1}}(0)=\C\1, \ W^{1,T_{1}}(\frac{1}{3})=\C y^1_{-\frac{1}{3}}\1=0, $$$$
 W^{1,T_{1}}(\frac{2}{3})=\C y^2_{-\frac{2}{3}}\1=\C y^2, \ W^{1,T_{1}}(\frac{4}{3})=\C y^1_{-\frac{4}{3}}\1=\C y^1,$$
 $$
 W^{2,T_{1}}(0)=\C w^2, \ W^{2,T_{1}}(\frac{1}{3})=\C y^1_{-\frac{1}{3}}w^2=\C w^1,
 $$$$
 W^{2,T_{1}}(\frac{2}{3})=\C  y^2_{-\frac{2}{3}}w^2=0, \ W^{2,T_{1}}(\frac{5}{3})=\C y^2_{-\frac{5}{3}}w^2=\C
 y^2(-2)w^2,
 $$
 $$
 L(0)|_{W^{1,T_{1}}(0)}=\frac{1}{36}id,  \
 L(0)|_{W^{2,T_{1}}(0)}=\frac{1}{9}id.
 $$
 \end{lem}
 We have the following result.
\begin{theorem}\label{twisted-module2}
$W^{1,T_{1}}$ and $W^{2,T_{1}}$ are the only two irreducible
$\sigma$-twisted modules of $V_{\Z\be}^{+}$.
\end{theorem}
\pf By Lemma \ref{twisted-module1}, $W^{1,T_{1}}$ and $W^{2,T_{1}}$
are inequivalent $\sigma$-twisted modules of $V_{\Z\be}^{+}$. Note
that $W^{1,T_{1}}$ and $W^{2,T_{1}}$ have irreducible quotients
which are $\sigma$-twisted modules of $V_{\Z\be}^{+}$ with lowest
weights $\dfrac{1}{36}$ and $\dfrac{1}{9},$ respectively. If
$W^{i,T_{1}}$ is not irreducible for some $i$, then the lowest
weight $\lambda$ of the maximal proper submodule is different from
$\dfrac{1}{36}$ and $\dfrac{1}{9}.$ By \cite{Y}, $V_{\Z\be}^{+}$ is
$C_2$-cofnite. It follows from \cite{DLM1} that $V_{\Z\be}^{+}$ has
an irreducible $\sigma$-twisted module with lowest weight $\lambda.$
This means that there are at least three inequivalent irreducible
$\sigma$-twisted modules of $V_{\Z\be}^{+}$, which contradicts Lemma
\ref{twisted-modules}. So $W^{1,T_{1}}$ and $W^{2,T_{1}}$ are
irreducible  inequivalent $\sigma$-twisted $V_{\Z\be}^{+}$-modules.
Then the theorem follows from Lemma \ref{twisted-modules}. \qed

Note that
$$
\sigma^{2}=\sigma^{-1}=e^{2\pi i(-h'(0))},
$$
and
$$
e^{2\pi i(-h'(0))}(-h')=-h', \ e^{2\pi
i(-h'(0))}y^1=\dfrac{-1-\sqrt{3}i}{2}y^1, \ e^{2\pi
i(-h'(0))}y^2=\dfrac{-1+\sqrt{3}i}{2}y^2.
$$
So we similarly have
\begin{lem}\label{twisted2}
$(W^{i,T_{2}},Y_{\sigma^{-1}}(\cdot,z))=(W^i,
Y(\Delta(-h',z)\cdot,z))$ are irreducible $\sigma^{-1}$-twisted
modules of $V_{L_{2}}, i=1,2.$
\end{lem}
 It is easy to see that
\begin{equation}\label{virasoro element2}
\Delta(-h',z)L(-2)\1=L(-2)\1-z^{-1}h'(-1)\1+\dfrac{1}{36}z^{-2}\1,
\end{equation}
\begin{equation}\label{s-element2}
Y_{\sigma^{-1}}(-h',z)=Y(-h'+\frac{1}{18}z^{-1},z),
\end{equation}
\begin{equation}\label{u-nilpotent2}
Y_{\sigma^{-1}}(y^1,z)=z^{-\frac{1}{3}}Y(y^1,z),
\end{equation}
\begin{equation}\label{l-nilpoten2}
Y_{\sigma^{-1}}(y^2,z)=z^{\frac{1}{3}}Y(y^2,z).
\end{equation}
By (\ref{virasoro element2})-(\ref{l-nilpoten2}), we have
$$ W^{1,T_{2}}(0)=\C\1, \ W^{1,T_{2}}(\frac{1}{3})=\C y^2_{-\frac{1}{3}}\1=0, $$$$
 W^{1,T_{2}}(\frac{2}{3})=\C y^1_{-\frac{2}{3}}\1=\C y^1, \ W^{1,T_{2}}(\frac{4}{3})=\C y^2_{-\frac{4}{3}}\1=\C y^2,$$
 $$
 W^{2,T_{2}}(0)=\C w^1, \ W^{2,T_{2}}(\frac{1}{3})=\C y^2_{-\frac{1}{3}}w^1=\C w^2,
 $$$$
 W^{2,T_{2}}(\frac{2}{3})=\C  y^1_{-\frac{2}{3}}w^1=0, \ W^{2,T_{2}}(\frac{5}{3})=\C y^1_{-\frac{5}{3}}w^1=\C
 y^1(-2)w^1,
 $$
 $$
 L(0)|_{W^{1,T_{2}}(0)}=\frac{1}{36}id,  \
 L(0)|_{W^{2,T_{2}}(0)}=\frac{1}{9}id.
 $$
Similar to Theorem \ref{twisted-module2}, we have
\begin{theorem}\label{twisted-module3}
$W^{1,T_{2}}$ and $W^{2,T_{2}}$ are the only two irreducible
$\sigma^{2}$-twisted modules of $V_{\Z\be}^{+}$.
\end{theorem}

We finally classify all the irreducible modules of
$V_{L_{2}}^{A_{4}}$. Recall that
$(V_{\Z\be}^{+})^{\<\sigma\>}=V_{L_{2}}^{A_{4}}$. We prove, in
particular, that any irreducible
$(V_{\Z\be}^{+})^{\<\sigma\>}$-module is contained in some
irreducible $\sigma^i$-twisted $V_{\Z\be}^{+}$-module, $i=0,1,2$.

Let $X^1$ and $X^2$ be defined as in (\ref{e3.9}). By Lemma
\ref{deco}, $X^i$ generates an irreducible
$(V_{\Z\be}^{+})^{\<\sigma\>}$-module with lowest weight 4, denoted
by $(V_{\Z\be}^{+})^{i}$, $i=1,2$.

Note that $W^{i,T_{1}}, W^{i,T_{2}}$, $i=1,2$ can also be regarded
as $(V_{\Z\be}^{+})^{\<\sigma\>}$-modules. Set
$$
w^{1,T_{1},1}=\1\in W^{1,T_{1}}(0), \ w^{1,T_{1},2}=y^2\in
W^{1,T_{1}}(\frac{2}{3}), \ w^{1,T_{1},3}=y^1\in
W^{1,T_{1}}(\frac{4}{3}),
$$
$$
w^{2,T_{1},1}=w^2\in W^{2,T_{1}}(0), \ w^{2,T_{1},2}=w^1\in
W^{2,T_{1}}(\frac{1}{3}), \ w^{2,T_{1},3}=y^2(-2)w^2\in
W^{2,T_{1}}(\frac{5}{3}),
$$
$$
w^{1,T_{2},1}=\1\in W^{1,T_{2}}(0), \ w^{1,T_{2},2}=y^1\in
W^{1,T_{2}}(\frac{2}{3}), \ w^{1,T_{2},3}=y^2\in
W^{1,T_{1}}(\frac{4}{3}),
$$
$$
w^{2,T_{2},1}=w^1\in W^{2,T_{2}}(0), \ w^{2,T_{2},2}=w^2\in
W^{2,T_{2}}(\frac{1}{3}), \ w^{2,T_{2},3}=y^1(-2)w^1\in
W^{2,T_{2}}(\frac{5}{3}).
$$

Then we have the following lemma.
\begin{lem}\label{irreducible modules3}
Let $W^{i,T_{j},k}$ be the  $(V_{\Z\be}^{+})^{\<\sigma\>}$-module
generated by $w^{i,T_{j},k}$, where $i,j=1,2$, $k=1,2,3$. Then
$W^{i,T_{j},k}$, $i,j=1,2,k=1,2,3$ are irreducible
$(V_{\Z\be}^{+})^{\<\sigma\>}$-modules such that
$$
L(0)w^{1,T_{1},1}=\frac{1}{36}w^{1,T_{1},1}, \
L(0)w^{1,T_{2},1}=\frac{1}{36}w^{1,T_{2},1},$$
$$
L(0)w^{1,T_{1},2}=\frac{25}{36}w^{1,T_{1},2}, \
L(0)w^{1,T_{2},2}=\frac{25}{36}w^{1,T_{2},2},
$$
$$L(0)w^{1,T_{1},3}=\frac{49}{36}w^{1,T_{1},3}, \
L(0)w^{1,T_{2},3}=\frac{49}{36}w^{1,T_{2},3},
$$
$$
L(0)w^{2,T_{1},1}=\frac{1}{9}w^{2,T_{1},1}, \
L(0)w^{2,T_{2},1}=\frac{1}{9}w^{2,T_{2},1},$$
$$
L(0)w^{2,T_{1},2}=\frac{4}{9}w^{2,T_{1},2}, \
L(0)w^{2,T_{2},2}=\frac{4}{9}w^{2,T_{2},2},
$$
$$L(0)w^{2,T_{1},3}=\frac{16}{9}w^{2,T_{1},3}, \
L(0)w^{2,T_{2},3}=\frac{16}{9}w^{2,T_{2},3}.
$$
\end{lem}
\pf The lemma follows from a general result: Let $U$ be a vertex
operator algebra with an automorphism $g$ of order $T.$ Let $M=\sum
_{n\in\frac{1}{T}\Z_+}M(n)$ be an irreducible $g$-twisted admissible
$U$-module. Then $M^i=\oplus_{n\in \frac{i}{T}+\Z}M(n)$ is an
irreducible $V^g$-module for $i=0,...,T-1$ (cf. \cite{DM1}). \qed

We have the following lemma from \cite{DM1}
\begin{lem}\label{decom1}
As an $(V_{\Z\be}^{+})^{\<\sigma\>}$-module,
$$V_{\Z+\frac{1}{4}\be}=(V_{\Z+\frac{1}{4}\be})^{0}\oplus
(V_{\Z+\frac{1}{4}\be})^{1}\oplus (V_{\Z+\frac{1}{4}\be})^{2}$$ such
that $(V_{\Z+\frac{1}{4}\be})^{0}$, $(V_{\Z+\frac{1}{4}\be})^{1}$
and $(V_{\Z+\frac{1}{4}\be})^{2}$ are irreducible
$(V_{\Z\be}^{+})^{\<\sigma\>}$-modules generated by
$e^{\be/4}+e^{-\be/4}$, $h(-2)\otimes
(e^{\be/4}-e^{-\be/4})-\sqrt{2}h(-1)^2\otimes
(e^{\be/4}+e^{-\be/4})+a(e^{3\be/4}+e^{-3\be/4})$ and $h(-2)\otimes
(e^{\be/4}-e^{-\be/4})-\sqrt{2}h(-1)^2\otimes
(e^{\be/4}+e^{-\be/4})-a(e^{3\be/4}+e^{-3\be/4})$ for some $0\neq
a\in\C$  with weights $\frac{1}{4}$, $\frac{9}{4}$ and $\frac{9}{4}$
respectively.
\end{lem}

We are now in a position to state the main result of this section.
Recall that $(V_{\Z\be}^{+})^{0}=(V_{\Z\be}^{+})^{\<\sigma\>}.$
\begin{theorem}
There are exactly  21 irreducible modules of
$(V_{\Z\be}^{+})^{\<\sigma\>}$. We give them by the following tables
1-4. \vskip 5ex
\begin{center}
\begin{tabular}{|c|c|c|c|c|c|c|}
\hline
&$(V_{\Z\be}^{+})^{0}$&$(V_{\Z\be}^{+})^{1}$&$(V_{\Z\be}^{+})^{2}$&
$V_{\Z\be}^{-}$&$V_{\Z\be+\frac{1}{8}\be}$&$V_{\Z\be+\frac{3}{8}\be}$\\
\hline
$\omega$&$0$&$4$&$4$&$1$&$\frac{1}{16}$&$\frac{9}{16}$\\
\hline
\end{tabular}
\end{center}

\begin{center}
\begin{tabular}{|c|c|c|c|c|c|c|}
\hline &$W^{1,T_{1},1}$&$W^{1,T_{1},2}$&$W^{1,T_{1},3}$&
$W^{2,T_{1},1}$&$W^{2,T_{1},2}$&$W^{2,T_{1},3}$\\
\hline
$\omega$&$\frac{1}{36}$&$\frac{25}{36}$&$\frac{49}{36}$&$\frac{1}{9}$&$\frac{4}{9}$&$\frac{16}{9}$\\
\hline
\end{tabular}
\end{center}

\begin{center}
\begin{tabular}{|c|c|c|c|c|c|c|}
\hline &$W^{1,T_{2},1}$&$W^{1,T_{2},2}$&$W^{1,T_{2},3}$&
$W^{2,T_{2},1}$&$W^{2,T_{2},2}$&$W^{2,T_{2},3}$\\
\hline
$\omega$&$\frac{1}{36}$&$\frac{25}{36}$&$\frac{49}{36}$&$\frac{1}{9}$&$\frac{4}{9}$&$\frac{16}{9}$\\
\hline
\end{tabular}
\end{center}

\begin{center}
\begin{tabular}{|c|c|c|c|}
\hline
&$(V_{\Z\be+\frac{1}{4}\be})^{0}$&$(V_{\Z\be+\frac{1}{4}\be})^{1}$&
$(V_{\Z\be+\frac{1}{4}\be})^{2}$\\
\hline $\omega$&$\frac{1}{4}$&$\frac{9}{4}$&$\frac{9}{4}$\\
\hline
\end{tabular}
\end{center}
\end{theorem}
\pf It follows from the proof of Lemma \ref{twisted-modules} and
Theorem 6.1 in \cite{DM1} that $V_{\Z\be}^{-}$,
$V_{\Z\be+\frac{1}{8}\be}$ and $V_{\Z\be+\frac{3}{8}\be}$ are
irreducible $(V_{\Z\be}^{+})^{\<\sigma\>}$-modules and as
$(V_{\Z\be}^{+})^{\<\sigma\>}$-modules,
$$ V_{\Z\be}^{-}\cong
V_{\Z\be+\frac{\be}{2}}^-\cong V_{\Z\be+\frac{\be}{2}}^+, \
V_{\Z\be+\frac{\be}{8}}\cong V_{\Z\be}^{T_{2},+}\cong
V_{\Z\be}^{T_{1},+}, \ V_{\Z\be+\frac{3\be}{8}}\cong
V_{\Z\be}^{T_{2},-}\cong V_{\Z\be}^{T_{1},-}.$$
 Then the theorem
follows from Lemma \ref{deco}, Lemma \ref{decom1}, Theorems
\ref{twisted-module2} and \ref{twisted-module3}, Theorem
\ref{rationality} and Theorem A in \cite{M1}. \qed

\end{document}